\documentclass[11pt]{article}

\usepackage{amssymb}
\usepackage{latexsym}

\newcommand{\bsq}{{\vrule height .9ex width .8ex depth -.1ex }}
\newcommand{\eeq}{\end{equation}}
\newcommand{\beql}[1]{\begin{equation}\label{#1}}
\newcommand{\beq}{\begin{equation}}
\newcommand{\eqn}[1]{(\ref{#1})}
\newcommand{\ZZ}{{\mathbb Z}}

\newtheorem{theorem}{Theorem}

\makeatletter
\def\@sect#1#2#3#4#5#6[#7]#8{\ifnum #2>\c@secnumdepth
     \def\@svsec{}\else
     \refstepcounter{#1}\edef\@svsec{\csname the#1\endcsname.\hskip .75em }\fi
     \@tempskipa #5\relax
      \ifdim \@tempskipa>\z@
        \begingroup #6\relax
          \@hangfrom{\hskip #3\relax\@svsec}{\interlinepenalty \@M #8\par}%
        \endgroup
       \csname #1mark\endcsname{#7}\addcontentsline
         {toc}{#1}{\ifnum #2>\c@secnumdepth \else
                      \protect\numberline{\csname the#1\endcsname}\fi
                    #7}\else
        \def\@svsechd{#6\hskip #3\@svsec #8\csname #1mark\endcsname
                      {#7}\addcontentsline
                           {toc}{#1}{\ifnum #2>\c@secnumdepth \else
                             \protect\numberline{\csname the#1\endcsname}\fi
                       #7}}\fi
     \@xsect{#5}}
\def\@begintheorem#1#2{\it \trivlist \item[\hskip \labelsep{\bf #1\ #2.}]}
\makeatother

\begin{document}

\begin{center}
{\Large\bf The Optimal Isodual Lattice Quantizer in Three Dimensions} \\

\vspace*{+.2in}
J. H. Conway\\
Mathematics Department \\
Princeton University \\
Princeton, NJ 08544 \\ [+.25in]

N. J. A. Sloane \\
AT\&T Shannon Labs \\
180 Park Avenue \\
Florham Park, NJ 07932-0971\\ 

\vspace*{+.1in}

Email: conway@math.princeton.edu,
njas@research.att.com \\
\vspace*{+.2in}
Jan 02 2006.
\vspace*{+.2in}

{\bf Abstract}
\vspace*{+.1in}
\end{center}

The mean-centered cuboidal (or m.c.c.) lattice is known to be the optimal packing and covering among all
isodual three-dimensional lattices.
In this note we show that it is also the best quantizer.
It thus joins the isodual lattices $\ZZ$, $A_2$ and (presumably) $D_4$,
$E_8$ and the Leech lattice in being simultaneously optimal with respect to all three criteria.

\vspace{0.8\baselineskip}

\vspace{0.8\baselineskip}

\setlength{\baselineskip}{1.5\baselineskip}
\section{Introduction}
An isodual lattice \cite{CoSl94} is one that is geometrically similar to its dual.
Let $\Lambda$ be an $n$-dimensional lattice and let $V$ denote the Voronoi cell containing the origin.
We may assume $\Lambda$ is scaled so that $V$ has unit volume.
Three important parameters of $\Lambda$ are the packing radius (the in-radius of $V$), the covering radius (the circum-radius of $V$) and its quantization
error, which is the normalized second moment
\beql{EqG}
G := \frac{1}{n} \int_V x \cdot x~ dx
\eeq
(see \cite{SPLAG}, also \cite{DW05}, \cite{GN98}, \cite{KB03}).

It is known that the face-centered cubic (or f.c.c.) lattice $A_3$ has the largest packing radius of any three-dimensional lattice (Gauss), while its dual,
the body-centered cubic (or b.c.c.) lattice $A_3^\ast$ has both the smallest
covering radius \cite{Bam54} and the smallest quantization error \cite{BS83}.
In \cite{CoSl94} it was shown that among all {\em isodual} three-dimensional
lattices, the mean-centered cuboidal (or m.c.c.) lattice
has the largest packing radius and the smallest covering radius.

The m.c.c. lattice, denoted here by $M_3$, has Gram matrix
$$\frac{1}{2} \left[
\begin{array}{ccc}
1+ \sqrt{2} & -1 & -1 \\ [+.1in]
-1 & 1+ \sqrt{2} & 1- \sqrt{2} \\ [+.1in]
-1 & 1- \sqrt{2} & 1+ \sqrt{2}
\end{array}
\right] \,,
$$
determinant 1, packing radius $\frac{1}{2} \sqrt{\frac{1}{2} + \frac{1}{\sqrt{2}}}$, center
density $\delta =0.1657 \ldots$ (which is 
between the values for the f.c.c. and b.c.c. lattices),
covering radius $3^{0.5} 2^{-1.25}$ (again between the values for the f.c.c. and b.c.c. lattices), kissing number 8 and automorphism group of order 16.
The m.c.c. lattice received a brief mention in \cite{BeM89} and was studied in more detail in \cite{CoSl94}.
It also arises from the period matrix of the hyperelliptic Riemann surface
$w^2 = z^8 -1$ \cite{BerSl97}.

The purpose of this note is to prove:
\begin{theorem}\label{th1}
The m.c.c. lattice $M_3$ is the optimal quantizer among all isodual
three-dimensional lattices.
\end{theorem}

In higher dimensions, less is known.
The lattice $D_4$ is optimal among all four-dimensional lattice packings, and has a lower quantization error than any other four-dimensional lattice presently known (see \cite{SPLAG} for references).
It is {\em not} the best four-dimensional lattice covering ($A_4^\ast$ is better), but it is isodual and may be optimal among isodual lattices with respect to all three criteria.
Similar remarks apply to the isodual eight-dimensional lattice $E_8$ (again
$A_8^\ast$ is a better covering but is not isodual).
In 24 dimensions the isodual Leech lattice is known to be the best lattice packing
\cite{CK04}, and may well also be the optimal covering and quantizer.

\section{Proof of Theorem \ref{th1}}
We will specify three-dimensional lattices by giving Gram matrices and conorms (or Selling parameters) --- cf. \cite{CoFu97},
\cite{CSLDL6}, \cite{CoSl94}.
It was shown in \cite{CoSl94} that, up to equivalence, indecomposable isodual three-dimensional lattices of determinant 1 have Gram matrices of the form
\beql{Eq1}
\frac{1}{2- \alpha \beta}
\left[
\begin{array}{ccc}
\frac{2 \alpha}{\beta} & - \alpha \beta & - \alpha (2- \beta ) \\ [.15in]
- \alpha \beta & \frac{2 \beta}{\alpha} & - \frac{2 \beta (1- \alpha )}{\alpha} \\ [+.1in]
- \alpha (2- \beta ) & \frac{- 2 \beta (1- \alpha )}{\alpha} & \frac{\alpha^2 \beta + 2 \alpha + 2 \beta - 4 \alpha \beta}{\alpha}
\end{array}
\right] \,,
\eeq
where $\alpha$, $\beta$ are any real numbers satisfying
$0 < \alpha < 1$, $0 < \beta < 1$; and decomposable lattices have Gram matrices
\beql{Eq2}
\left[ \begin{array}{ccc}
1 & 0 & 0 \\
0 & \alpha & -h \\
0 & - h & \beta
\end{array}
\right] \,,
\eeq
where $\alpha$, $\beta$, $h$ are any real numbers satisfying $0 \le 2h \le \alpha \le \beta$, $\alpha \beta - h^2 =1$.

In the indecomposable case the nonzero conorms are:
\begin{eqnarray}\label{Eq3}
p_{01} & = & \frac{\alpha (2- \beta )}{\gamma } , ~~
p_{02} = \frac{\alpha \beta}{\gamma } , ~~ p_{03} = \frac{2 \alpha (1- \beta )}{\beta \gamma } \,, \nonumber \\ [+.15in]
p_{12} & = & \frac{2 \beta (1- \alpha )}{\alpha \gamma } , ~~
p_{13} = \frac{2 (1- \alpha ) (1- \beta )}{\gamma } , ~~
p_{23} = \frac{\beta (2- \alpha )}{\gamma } \,,
\end{eqnarray}
where $\gamma  = 2 - \alpha \beta$; in the decomposable case they are:
\beql{Eq4}
p_{01}  =  1, ~~ p_{02} = \alpha -h , ~~ p_{03} = \beta -h , ~~
p_{12}  =  0 , ~~ p_{13} =0 , ~~p_{23} = h \,.
\eeq

The m.c.c. lattice corresponds to the case $\alpha = \beta = 2 - \sqrt{2}$ of Eqs. \eqn{Eq1} and \eqn{Eq3}.

From \cite{BS83} we know that the normalized second moment \eqn{EqG} for a decomposable or indecomposable three-dimensional lattice $\Lambda$ is given by
\beql{Eq5}
G = \frac{DS_1 + 2S_2 +K}{36 D^{4/3}} \,,
\eeq
where $D= \det \, \Lambda$,
\begin{eqnarray*}
S_1 & = & p_{01} + p_{02} + p_{03} + p_{12} + p_{13} + p_{23} \,, \\
S_2 & = & p_{01} p_{02} p_{13} p_{23} + p_{01} p_{03} p_{12} p_{23} + p_{02}
p_{03} p_{12} p_{13} \,, \\
K & = &
p_{01} p_{02} p_{03} (p_{12} + p_{13} + p_{23} )
+ p_{01} p_{12} p_{13} (p_{02} + p_{03} + p_{23} ) \\
&& + p_{02} p_{12} p_{23} (p_{01} + p_{03} + p_{13} ) 
+ p_{03} p_{13} p_{23} (p_{01} + p_{02} + p_{12} ) \,.
\end{eqnarray*}
For our lattices, $D=1$.

We first consider the indecomposable case.
From \eqn{Eq3}, \eqn{Eq5} we find that
\beql{Eq6}
G = \frac{f ( \alpha , \beta )}{36 \alpha \beta (2- \alpha \beta )^4}
\eeq
where
$$\begin{array}{l@{\,}l}
& f(\alpha , \beta ) = 3 \alpha^5 \beta^5 - 8 \alpha^4 \beta^4 (\alpha + \beta ) +
4 \alpha^3 \beta^3 (\alpha^2 + 10 \alpha \beta + \beta^2 ) \\ [+.05in]
- & 48 \alpha^3 \beta^3 ( \alpha + \beta ) + 8 \alpha^2 \beta^2 (3
\alpha^2 +4 \alpha \beta + 3 \beta^2 ) + 32 \alpha^2 \beta^2 ( \alpha + \beta ) \\ [+.05in]
- & 8 \alpha \beta (5 \alpha^2 + 6 \alpha \beta + 5 \beta^2 ) + 16
(\alpha^2 + \alpha \beta + \beta^2 ) \,.
\end{array}
$$
It turns out that there is exactly one choice for $(\alpha , \beta )$ in the range $0 < \alpha < 1$, $0 <\beta < 1$ for which both partial derivatives
$\frac{\partial G}{\partial \alpha}$ and $\frac{\partial G}{\partial \beta}$ vanish,
namely $\alpha = \beta = 2 - \sqrt{2}$, corresponding to the m.c.c. lattice.
At all other points in the interior of this region, one of the two partial
derivatives does not vanish and so the point cannot be a local minimum.
To show this we used the computer algebra system Maple \cite{Map1}
to compute the numerators of
$\frac{\partial G}{\partial \alpha}$ and $\frac{\partial G}{\partial \beta}$, giving
a pair of simultaneous equations in $\alpha$ and $\beta$, symmetrical in $\alpha$ and $\beta$.
By eliminating $\beta$, we obtain a single equation for $\alpha$:
\begin{eqnarray*}
&& \alpha ( \alpha -1 ) ( \alpha -2 ) ( \alpha^2 +2 \alpha -2 ) ( \alpha^2 -4 \alpha +2 ) 
 (\alpha^4 - 4 \alpha^3 + 6 \alpha^2 -4 ) ~ \times  \\
&& \times ~ (57 \alpha^6 - 220 \alpha^5 - 102 \alpha^4 + 1448 \alpha^3 - 1860 \alpha^2 + 832 \alpha - 152 ) =0 \,.
\end{eqnarray*}
The only roots in the range $0 < \alpha < 1$ are $2- \sqrt{2}$, $\sqrt{3} -1$ and the root $0.9894 \ldots$ of the sixth-degree factor.
By symmetry, $\beta$ must also take one of these three values.
At only one of these nine points do both partial derivatives vanish, namely
$\alpha = \beta = 2- \sqrt{2}$.
At that point the matrix of second partial derivatives is positive definite, showing that m.c.c. is a local minimum.
The resulting value of $G$ is
$$
\frac{17 + 4 \sqrt{2}}{288}
= 0.0786696 \ldots \,,$$
between the values for the b.c.c. and f.c.c. lattices, which are respectively
$$\frac{19}{384} 2^{2/3} = 0.0785432 \ldots ~~\mbox{and} ~~
\frac{2^{1/3}}{16} = 0.0787450 \ldots \,.
$$

On the boundary of the region the lattices are either degenerate
(if $\alpha$ or $\beta$ is 0) or decomposable (if $\alpha$ or $\beta$ is $1$).
In the latter case we assume $\alpha \le \beta$ and find that there is a unique
point where $\frac{\partial G}{\partial \alpha}$ 
and $\frac{\partial G}{\partial \beta}$ vanish,
when $\alpha = \sqrt{3} -1$, $\beta =1$.
This is the lattice $\ZZ \oplus 3^{-1/4} A_2$, for which $G= \frac{5 \sqrt{3}}{162} + \frac{1}{36} = 0.0812361 \ldots$.
It is not a local minimum.

It remains to consider the decomposable case.
From \eqn{Eq4}, \eqn{Eq5}, we find that
\beql{Eq76}
G = \frac{1}{36}
\left\{ \alpha \beta ( \alpha + \beta ) + 
2(\alpha \beta -1 )^{3/2} + 2 \alpha + 2 \beta +1 \right\} \,.
\eeq
Again we solve $\frac{\partial G}{\partial \alpha} = \frac{\partial G}{\partial \beta} =0$, and
find that the only possibilities in the range $0 \le \alpha \le \beta$ are
$\alpha = \beta =1$, $G= \frac{1}{12}$;
$\alpha =1$, $\beta =2$, $G= \frac{1}{12}$; and
$\alpha = \beta = \sqrt{3} -1$ (the decomposable lattice mentioned above).
None of these are local minima.
This completes the proof.
The proof also shows that the m.c.c. lattice is the only isodual lattice where $G$ has a local minimum.~~~$\bsq$


\end{document}